# On Jacobi's condition for the simplest problem of calculus of variations with mixed boundary conditions

Milan Batista

University of Ljubljana, Faculty of Maritime Studies and Transport, Slovenia

milan.batista@fpp.uni-lj.si

(June 2015)

**Abstract**

The purpose of this paper is an extension of Jacobi's criteria for positive definiteness of second variation of the simplest problems of calculus of variations subject to mixed boundary conditions. Both non constrained and isoperimetric problems are discussed. The main result is that if we stipulate conditions (21) and (22) then Jacobi's condition remains valid also for the mixed boundary conditions.

**1 Introduction**

As it is well known [1; 2; 3; 4; 5; 6; 7] the simplest problem of the calculus of variations is to find a function $y = y(x)$ that minimize (or maximize) the functional

$$J[y] \equiv \int_a^b F(x, y, y') dx \qquad (1)$$

under Dirichlet boundary conditions

$$y(a) = A \text{ (fixed)}, \quad y(b) = B \text{ (fixed)} \qquad (2)$$

where function $F$ and constants $A$ and $B$ are given, $(\ )' \equiv d(\ )/dx$. For convenience we shall call this problem the Dirichlet problem.

**Note.** In the paper we shall assume that for all the functions we are going to use the domain of definition is the interval $[a,b]$ and that they possess continuous derivatives with respect to all its arguments as many order as needed, unless stated otherwise. We introduce a function by writing $y = y(x)$ (as example) and in the sequel use $y$ to denote the function and $y(x)$ to denote its value.

Suppose that $y$ is an extremal of $J[y]$. Then a sufficient condition for $y$ to realize a weak minimum of $J[y]$ is that its second variation $\delta^2 J[h]$ is positive definite, that is, $\delta^2 J[h] > 0$ for any piecewise smooth, i.e., continuous and piecewise continuously differentiable, $h = h(x) \not\equiv 0$ such that

$$h(a)=0, \quad h(b)=0 \tag{3}$$

In this paper we will assume that $\delta^2 J[h]$ is given in the integrated form

$$\delta^2 J[h] \equiv Rh^2 \Big|_a^b + \int_a^b \left( Ph'^2 + Qh^2 \right) dx \tag{4}$$

where

$$P(x) \equiv \frac{\partial^2 F}{\partial y'^2}, \quad Q(x) \equiv \frac{\partial^2 F}{\partial y^2} - \frac{dR}{dx}, \quad R(x) \equiv \frac{\partial^2 F}{\partial y \partial y'}, \tag{5}$$

Here and in the sequel the partial derivatives are evaluated at $(x,y,y')$. We abbreviate the condition $\delta^2 J[h]>0$ as $\delta^2 J>0$ to mean that it must be fulfilled for piecewise smooth function $h=h(x)\not\equiv 0$ which satisfy the constraints of the given variation problem. We will call such $h$ an admissible variation.

Now the question under what conditions we have $\delta^2 J>0$, the answer is the following Jacobi's theorem (see for instance [1; 4; 5; 6])

**Theorem 1 (Dirichlet problem).** *In order that $\delta^2 J>0$ it is necessary and sufficient that the following conditions hold:*

1. $P(x)>0$ for all $a \leq x \leq b$ (strengthened Legendre condition)
2. $u(x) \neq 0$ for all $a < x \leq b$ (strengthened Jacobi condition)

*where $u=u(x)$ is the solution of Jacobi's accessory equation*

$$\mathbf{L}(u) \equiv -\frac{d}{dx}(Pu') + Qu = 0 \tag{6}$$

*which satisfy the initial conditions*

$$u(a)=0, \quad u'(a)=1 \tag{7}$$

Note that $u'(a)=1$ is only for sake of definiteness of $u$ [5].

For a simplest isoperimetric problem $y$ must beside boundary conditions (2) satisfy also the additional condition

$$\int_a^b G(x,y,y') dx = \ell \tag{8}$$

where $\ell$ is a given constant. In this case an extremal of the problem $y$ is obtained from the functional

$$J[y] = \int_a^b H(x,y,y') dx \tag{9}$$

where $H \equiv F + \lambda G$, $\lambda$ is a Lagrange multiplier. The requirement for $\delta^2 J>0$ which is given by (4) and where (5) is replaced by

$$P(x) \equiv \frac{\partial^2 H}{\partial y'^2}, \qquad Q(x) \equiv \frac{\partial^2 H}{\partial y^2} - \frac{dR}{dx}, \qquad R(x) \equiv \frac{\partial^2 H}{\partial y \partial y'}. \tag{10}$$

is in this case enhanced by requirement that an admissible $h$ must also satisfy the condition

$$\left. \frac{\partial G}{\partial y'} h \right|_a^b + \int_a^b hT\, dx = 0 \tag{11}$$

where

$$T(x) \equiv \frac{\partial G}{\partial y} - \frac{d}{dx} \frac{\partial G}{\partial y'} \tag{12}$$

For isoperimetric problem the answer to the question for the conditions for $\delta^2 J > 0$ is the following Bolza's theorem [8]

**Theorem 2 (Dirichlet isoperimetric problem).** *Let $P(x) > 0$ and $T \not\equiv 0$ on $[a,b]$. Let $u = u(x)$ and $v = v(x)$ be two solutions of*

$$\mathbf{L}(u) = 0, \qquad \mathbf{L}(v) = T \tag{13}$$

*respectively, that satisfy*

$$u(a) = 0, \qquad v(a) = 0 \tag{14}$$

*Further let*

$$m(x,a) \equiv \int_a^x uT\, dx, \qquad n(x,a) \equiv \int_a^x vT\, dx \tag{15}$$

*and let*

$$\Delta(x,a) \equiv mv - nu \tag{16}$$

*Then the necessary and sufficient condition for $\delta^2 J > 0$ is that $\Delta(x,a) \neq 0$ for all $a < x \leq b$.*

Absence of proof of Jacobi criteria for the mixed boundary conditions (what we shall call the mixed problem) of the form

$$y(a) = A \text{ (fixed)}, \qquad y'(b) = 0 \text{ (free)} \tag{17}$$

in textbooks was noted by Maddocks [9]. He state that the Jacobi's condition for $\delta^2 J > 0$ is valid also for the unconstraint mixed problem for all admissible $h$ such that

$$h(a) = 0, \qquad h'(b) = 0 \tag{18}$$

He note that the proof is slight generalization of the standard account of conjugate points given in [5]. The mixed problem Maddocks consider was the problem of stability of elastic rods. For more on mixed problem in stability of rods see for instance [10; 11; 12] and for broader aspect [13; 14].

In this paper we will establish a necessary and sufficient condition for $\delta^2 J > 0$ when extremal $y = y(x)$ of $J[y]$ satisfy the mixed boundary conditions

$$y'(a) = 0 \text{ (free)}, \quad y(b) = B \text{ (fixed)} \tag{19}$$

In this case an admissible variation $h$ must satisfy the following boundary conditions

$$h'(a) = 0, \quad h(b) = 0 \tag{20}$$

We will consider both unconstrained and isoperimetric case of this mixed problem. The proofs given below are chiefly based on work of Bolza [8; 1].

## 2 Mixed problem

Before we turn to the theorems, we consider the boundary terms $Rh^2\big|_a^b$ in (4) and $\frac{\partial G}{\partial y'} h \big|_a^b$ in (11).

These terms vanish identically for Dirichlet boundary conditions (3). However, for mixed boundary condition it vanish only at end where $h$ is zero. In order to take these terms completely vanish for mixed boundary conditions (20) it is thus sufficient to assume that

$$R(a) = 0 \tag{21}$$

$$\frac{\partial G}{\partial y'}(a) = 0 \text{ (isoperimetric problem only)} \tag{22}$$

If (21) and (22) holds then (4) and (11) becomes

$$\delta^2 J = \int_a^b \left( Ph'^2 + Qh^2 \right) dx \tag{23}$$

$$\int_a^b hT \, dx = 0 \tag{24}$$

We will restrict our discussion to the mixed problems where (21) and (22) holds and therefore (23) and (24) are valid.

In what follows we will use the following lemma

**Lemma 1.** *If h is admissible variation for the mixed problem satisfying* (20) *and if* (21) *holds then*

$$\delta^2 J = \int_a^b h \mathbf{L}(h) dx \tag{25}$$

Proof. Integration by parts (23) give

$$\delta^2 J = Phh'\big|_a^b + \int_a^b h \mathbf{L}(h) dx \tag{26}$$

By (20) boundary term vanish so we are left with (25).

2.1 Unconstrained mixed problem

To prove Theorem 3 we need the following Jacobi's lemma

**Lemma 2.** If $u = u(x)$ is a solution of (6) and $p = p(x)$ is an arbitrary pricewise smooth function, then

$$P(pu)'^2 + Q(pu)^2 = P(p'u)^2 + \frac{d}{dx}(Pp^2 uu') \qquad (27)$$

The proof can be found in Bolza [1]

The answer to the question for conditions that $\delta^2 J > 0$ for the simplest variation problem with mixed boundary (19) is the following:

**Theorem 3 (mixed problem).** Let $P(x) > 0$ on $[a,b]$ and let (21) holds. Let $u = u(x)$ be nontrivial solutions of

$$\mathbf{L}(u) = 0 \qquad (28)$$

that satisfy

$$u'(a) = 0 \qquad (29)$$

Then the necessary and sufficient condition for $\delta^2 J > 0$ is that $u(x) \neq 0$ for all $a < x \leq b$.

A definiteness of $u$ is explained in Appendix.

*Proof.* That the condition is necessary can be proven by opposite assertion: if $u(c) = 0$ for some $a < c \leq b$ then we can found $h$ that makes $\delta^2 J \leq 0$. Consider piecewise smooth function

$$h \equiv \begin{cases} u & a \leq x \leq c \\ 0 & c < x \leq b \end{cases} \qquad (30)$$

Clearly $h(b) = 0$ and since $u'(a) = 0$ also $h'(a) = 0$. Thus this $h$ is admissible variation and therefore obviously fulfil the condition of Lemma 1. Substituting (30) into (25) gives

$$\delta^2 J = \int_a^c u\mathbf{L}(u)dx \qquad (31)$$

But $\mathbf{L}(u) = 0$ and consequently $\delta^2 J = 0$.

To proof that the condition is sufficient we start with assumption $u(x) \neq 0$ for all $a \leq x \leq b$. With an admissible $h$ we can then define $p = p(x)$ by

$$p \equiv \frac{h}{u} \qquad (32)$$

Since $h$ is piecewise smooth, so is $p$. Using this $p$ in identity (27) and integrating both sides from $a$ to $b$ we find that

$$\delta^2 J = \int_a^b \left( Ph'^2 + Qh^2 \right) dx = Pp^2 uu' \Big|_a^b + \int_a^b P(p'u)^2 dx \tag{33}$$

Since $h(b)=0$ also $p(b)=0$. From $p'=(h'u-hu')/u^2$ and $h'(a)=u'(a)=0$ we get $p'(a)=0$. Therefore the boundary term vanish. If $P>0$, then $\delta^2 J>0$ because $(p'u)^2 > 0$ unless $p'u \equiv 0$. However, in that case $p'=0$ which together with $p(b)=0$ imply $p \equiv 0$. But this would be possible only if $h \equiv 0$ (no variation), contrary to assumption that $h$ is an admissible variation. This completes the proof.

2.2 Isoperimetric mixed problem

The following variant of Bolza's lemma [8] is needed to prove Theorem 4.

**Lemma 3.** *Let $u=u(x)$ and $v=v(x)$ be a solutions of*

$$\mathbf{L}(u)=0, \quad \mathbf{L}(v)=T \tag{34}$$

*respectively, that satisfy*

$$u'(a)=0, \quad v'(a)=0 \tag{35}$$

*and let $p=p(x)$ and $q=q(x)$ be two an arbitrary piecewise smooth function. Then*

$$P(pu+qv)'^2 + Q(pu+qv)^2 = P(p'u+q'v)^2 - 2q(p'm+q'n) + \frac{d}{dx}\left[ P(pu+qv)(pu'+qv') + (pm+qn)q \right] \tag{36}$$

where $m$ and $n$ are given by (15).

*Proof.* By direct calculation we obtain, using (34), (see [8] for details)

$$P(pu+qv)'^2 + Q(pu+qv)^2 = P(p'u+q'v)^2 - P(p'q-pq')(u'v-uv') + (pu+qv)qT + \frac{d}{dx}\left[ P(pu+qv)(pu'+qv') \right] \tag{37}$$

Two additional identities are needed to deduce (36) from this. The first one is

$$(pm'+qn')q = \frac{d}{dx}\left[ (pm+qn)q \right] - (p'm+q'n)q - (pm+qn)q' \tag{38}$$

The second one is obtained from the identity $v\mathbf{L}(u)-u\mathbf{L}(v)=-uT$ or explicitly

$$\frac{d}{dx}\left[ P(u'v-uv') \right] = uT \tag{39}$$

Integrating from $a$ to $x$ we obtain, using (15),

$$P(u'v-uv') = P(u'v-uv')\big|_{x=a} + m \tag{40}$$

By (35) the boundary term vanish so we are left with

$$P(u'v - uv') = m \tag{41}$$

Combining (37), (41) and (38) gives (36).

We note that Bolza prove this identity assuming $u(a) = 0$ and $v(a) = 0$. However, we see from (40) that the identity holds for any combinations of initial conditions.

For isoperimetric mixed problem the answer to the question for the conditions for $\delta^2 J > 0$ is the now the following modification of the Bolza's theorem

**Theorem 4 (mixed isoperimetric problem).** *Let $P(x) > 0$ and $T(x) \not\equiv 0$ on $[a,b]$ and let (21) and (22) holds. Let $u = u(x)$ and $v = v(x)$ be nontrivial solutions of*

$$\mathbf{L}(u) = 0, \quad \mathbf{L}(v) = T \tag{42}$$

*respectively, that satisfy*

$$u'(a) = 0, \quad v'(a) = 0 \tag{43}$$

*Further let*

$$\Delta(x,a) \equiv mv - nu \tag{44}$$

where $m$ and $n$ are given by (15). *Then the necessary and sufficient condition for $\delta^2 J > 0$ is that $\Delta(x,a) \neq 0$ for all $a < x \leq b$.*

A construction of $u$ and $v$ from the general solution of nonhomogeneous Jacobi's accessory equation is presented in Appendix.

*Proof.* Again, to show that the condition is necessary we suppose opposite, that is, if $\Delta(c,a) = 0$ for some $a < c \leq b$ then we can setup $h$ that makes $\delta^2 J \leq 0$. Define piecewise smooth function

$$h \equiv \begin{cases} \alpha v - \beta u & a \leq x \leq c \\ 0 & c < x \leq b \end{cases} \tag{45}$$

where

$$\alpha \equiv m(c,a), \quad \beta \equiv n(c,a) \tag{46}$$

From (45) we have $h(b) = 0$ and also, together with (43), $h'(a) = 0$. Moreover

$$\int_a^b hT\,dx = \alpha \int_a^c vT\,dx - \beta \int_a^c uT\,dx = \alpha n(c,a) - \beta m(c,a) = 0 \tag{47}$$

Conditions (11) and (20) are thus satisfied so this $h$ is an admissible variation. Substituting (45) into (25), we have

$$\delta^2 J = \int_a^b h\big[\alpha \mathbf{L}(v) - \beta \mathbf{L}(u)\big]dx \tag{48}$$

By (42) this reduces to

$$\delta^2 J = \alpha \int_a^b hT\,dx \tag{49}$$

But by (47) this integral vanish, hence $\delta^2 J = 0$.

To proof that the condition is also sufficient we start by assumption $\Delta(x,a) \neq 0$ for all $a \leq x \leq b$. Following Bolza [8] we define $p = p(x)$ and $q = q(x)$ by

$$pu + qv = h \tag{50}$$

$$mp' + nq' = 0 \tag{51}$$

where $h$ is an admissible variation ( note that (51) eliminates non-quadratic term on the right side of identity (36)) . The identity

$$\frac{d}{dx}(pm + qn) = (p'm + q'm) + (pm' + qn') \tag{52}$$

on using (50), (51) and (15) become

$$\frac{d}{dx}(pm + qn) = hT \tag{53}$$

Integrating subject to $m(a,a) = n(a,a) = 0$ yields

$$pm + qn = \int_a^x hT\,dx \tag{54}$$

From (50) and (54) we now obtain

$$p = \frac{v \int_a^x hT\,dx - nh}{\Delta}, \quad q = \frac{mh - u \int_a^x hT\,dx}{\Delta} \tag{55}$$

Both these functions are pricewise smooth on $(a,b]$, because $h$ is. At $x = a$ where $\Delta(a,a) = 0$ both becomes $\frac{0}{0}$. We will evaluate this indeterminate form by L'Hopital's rule. The first nonzero derivative of $\Delta$ at $x = a$ is, assuming $T(a) \neq 0$,

$$\Delta'''(a,a) = -2\frac{u(a)T^2(a)}{P(a)} \tag{56}$$

This derivative was evaluated by means of (44) and (41). The third derivatives of numerators of (55) at $x = a$ are

$$\left(v\int_a^x hT\,dx - nh\right)'''\bigg|_{x=a} = 2(hv'' - h''v)T\big|_{x=a}, \quad \left(mh - u\int_a^x hT\,dx\right)'''\bigg|_{x=a} = 2(h''u - hu'')T\big|_{x=a} \tag{57}$$

Using these and (56) we obtain

$$p(a) = -\frac{(hv'' - h''v)P}{uT}\bigg|_{x=a}, \quad q(a) = \frac{(hu'' - h''u)P}{vT}\bigg|_{x=a} \tag{58}$$

Now, using these $p$ and $q$ in (36) and integrate both sides from $a$ to $b$ we get

$$\delta^2 J = \int_a^b \left(Ph'^2 + Qh^2\right)dx = \left[Ph(pu' + qv') + q\int_a^x hT\,dx\right]_a^b + \int_a^b P(p'u + q'v)^2\,dx \tag{59}$$

The boundary term of this expression is zero because $u'(a) = v'(a) = 0$ and $\int_a^a hT\,dx = 0$, and because $h(b) = 0$ and $\int_a^b hT\,dx = 0$. If $P > 0$ then $\delta^2 J \geq 0$ where $\delta^2 J = 0$ only if $p'u + q'v \equiv 0$. However, that, together with (51), imply $p = \alpha$ and $q = \beta$ where $\alpha, \beta$ are constants. From (55) when $x = b$ we get $\alpha = \beta = 0$, thus $p = q = 0$. But when $\Delta \neq 0$ the trivial solution of (50) and (54) is possible only if $h \equiv 0$, which contradicts assumption that $h$ is admissible variation. Therefore $\delta^2 J > 0$, which is the desired conclusion.

## 3 Examples

In this section we will give two examples of use of Theorem 3 and Theorem 4 in the case of catenary where we have to find extreme of the integral (potential energy)

$$J[y] = \int_a^b y\sqrt{1 + y'^2}\,dx \tag{60}$$

subject to boundary condition

$$y'(a) = 0, \quad y(b) = y_b \tag{61}$$

Here $F = y\sqrt{1 + y'^2}$ and the extremal is

$$y = \omega \cosh\frac{x - a}{\omega} \tag{62}$$

where $\omega$ is a solution of

$$y_b = \omega \cosh\frac{b - a}{\omega} \tag{63}$$

To find which solution of (63) realize minimum of (60) we use Theorem 3. By straightforward computations (5) becomes

$$P = \frac{\omega}{\cosh^2\left(\frac{x-a}{\omega}\right)}, \quad Q = -\frac{1}{\omega\cosh^2\left(\frac{x-a}{\omega}\right)}, \quad R = \tanh\left(\frac{x-a}{\omega}\right). \tag{64}$$

Hence $R(a) = 0$ while we must choose $\omega > 0$ to ensure $P > 0$. By (64) Jacobi's accessory equation (6) become

$$\frac{d^2u}{dx^2} + \frac{2}{\omega}\tanh\left(\frac{x-a}{\omega}\right)\frac{du}{dx} + \frac{u}{\omega^2} = 0. \tag{65}$$

Solution of this equation satisfying $u'(a) = 0$ is

$$u = \left(\frac{x-a}{\omega}\right)\sinh\left(\frac{x-a}{\omega}\right) - \cosh\left(\frac{x-a}{\omega}\right). \tag{66}$$

where other constant is set to one. Equation $u(\zeta) = 0$ where $\zeta \equiv \frac{x-a}{\omega}$, has solution $\zeta \approx 1.1998$.

The condition $\delta^2 J > 0$ will be, according to Theorem 3, fulfil if $b < a + \zeta\omega$, i.e., if $\omega > \frac{b-a}{\zeta}$. Such $\omega$ is the largest of two possible solutions of (63). In other words, between two possible catenaries satisfying (61) the one with higher apex (see (62)) has minimum of potential energy, as in the case of Dirichlet conditions [4].

Consider now a catenary of a given length $\ell$. We have, in addition to (60) and (61),

$$\int_a^b \sqrt{1+y'^2}\, dx = \ell \tag{67}$$

so $G = \sqrt{1+y'^2}$ and $H = y\sqrt{1+y'^2} + \lambda\sqrt{1+y'^2}$. The extremal is

$$y = y_b - \omega\left[\cosh\left(\frac{b-a}{\omega}\right) - \cosh\left(\frac{x-a}{\omega}\right)\right] \tag{68}$$

where $\omega$ is a solution of

$$\omega\sinh\left(\frac{b-a}{\omega}\right) = \ell \tag{69}$$

This equation has only one solution for $\omega > 0$. To see if this solution give minimum of (60) we use Theorem 4. Using $H$ we from (10) for $P$, $Q$, $R$ obtain expressions given by (64). From (12) we have

$$\frac{\partial G}{\partial y'} = \tanh\left(\frac{x-a}{\omega}\right), \quad T = -\frac{1}{\omega\cosh^2\left(\frac{x-a}{\omega}\right)} \tag{70}$$

Hence $\frac{\partial G}{\partial y'}(a) = 0$ and $T(a) = -\frac{1}{\omega} \neq 0$, since $\omega > 0$ by assumption. With $T$ given by (70) the solution of nonhomogeneous Jacobi's accessory equation is

$$v = 1 + \left(\frac{x-a}{\omega}\right)\sinh\left(\frac{x-a}{\omega}\right) - \cosh\left(\frac{x-a}{\omega}\right). \tag{71}$$

By means of (66), (71) and (70) the expression (44) becomes

$$\Delta = \frac{(\cosh\zeta - \zeta\sinh\zeta)(\cosh(2\zeta) - \sinh(2\zeta)) - \zeta\sinh\zeta - \cosh\zeta(1-2\zeta)}{\cosh(2\zeta) - \sinh(2\zeta) + 1}. \tag{72}$$

where $\zeta \equiv \dfrac{x-a}{\omega}$. But only solution of $\Delta(\zeta)=0$ is $\zeta=0$ imply zero at $x=a$. Hence, by Theorem 4, we have $\delta^2 J > 0$, i.e. the minimum of potential energy, as expected.

### 4 Conclusion

It was shown that Jacoby's condition for positives of second variation of the simplest vibrational problem both unconstrained or isoperimetric remains valid when conditions (21) and (22) holds. In this respect, we note that (21) is automatically satisfied for a problems where $\dfrac{\partial^2 F}{\partial y \partial y'} \equiv 0$, i.e., when $F = f_0(x, y') + f_1(x, y)$ (like in the case of elastic rods), or, as we demonstrate in the case of catenary, when $R$ is an odd function with respect to $x = a$. Similarly for isoperimetric problems (22) is automatically satisfied when $G = G(x, y)$, or when $\dfrac{\partial G}{\partial y'}$ is odd function regard $x = a$.

### Appendix.

In this appendix we will show how $u$ and $v$, used in the Theorem 4, can be obtained from the general solution of nonhomogeneous Jacobi's accessory equation

$$\mathbf{L}(w) = \mu T \qquad (73)$$

where $\mu$ is a constant. This is second order linear differential equation which has the general solution in the form [15]

$$w = \mu \theta_0 + c_1 \theta_1 + c_2 \theta_2 \qquad (74)$$

Here $c_1$ and $c_2$ are arbitrary constants, $\theta_0 = \theta_0(x)$ is solution of $\mathbf{L}(\theta_0) = T$, and $\theta_1 = \theta_1(x)$, $\theta_2 = \theta_2(x)$ are linearly independent solutions of homogeneous part of (73). The solution (74) can equivalently be expressed by

$$w = \nu u + \mu v \qquad (75)$$

where $\nu \neq 1$ is a constant, $u = u(x)$ is solution of $\mathbf{L}(u) = 0$ and $v = v(x)$ is solution of $\mathbf{L}(v) = T$. To express $u$ and $v$ by $\theta_0, \theta_1$ and $\theta_2$ we set

$$c_1 = \nu \overline{C}_1 + \mu C_1 \qquad c_2 = \nu \overline{C}_2 + \mu C_2 \qquad (76)$$

where $\overline{C}_1, \overline{C}_2, C_1, C_2$ are constants. Combining (74) and (75) we obtain, by equate the terms in $\mu$ and $\nu$,

$$u = \overline{C}_1 \theta_1 + \overline{C}_2 \theta_2 \qquad (77)$$

$$v = \theta_0 + C_1 \theta_1 + C_2 \theta_2 \qquad (78)$$

These functions must satisfy initial conditions $u'(a)=0$ and $v'(a)=0$. Thus both of them may contain a constant. Clearly, such constant cannot affect calculation of $\delta^2 J$ but for the sake of definiteness it is desired to be specified in some way. One way is by normalization conditions $u(a)=1$ and $v(a)=1$. This is always possible because $\theta_1$ and $\theta_2$ are linearly independent solution of (73). Another way, suggested by Bolza [1], is the following. From $u'(a)=0$ we have $\overline{C}_1\theta_1'(a)+\overline{C}_2\theta_2'(a)=0$ and this is satisfied identically with $\overline{C}_1=\theta_2'(a)$ and $\overline{C}_2=-\theta_1'(a)$, hence

$$u=\theta_2'(a)\theta_1(x)-\theta_1'(a)\theta_2(x) \tag{79}$$

The constants $C_1$ and $C_2$ must be chosen in a way to satisfy condition

$$C_1\theta_1'(a)+C_2\theta_2'(a)+\theta_3'(a)=0 \tag{80}$$

The same methods can be used for *u* appearing in the Theorem 3.